\documentclass[a4paper,12pt,reqno]{amsart}
\usepackage{amssymb, amscd, epsfig}

\newtheorem{theor}{Theorem}
\newtheorem{prop}[theor]{Proposition}
\newtheorem{lem}[theor]{Lemma}
\theoremstyle{definition}
\newtheorem{cor}[theor]{Corollary}
\theoremstyle{remark}
\newtheorem{rem}{Remark}
\newtheorem{exmp}{Example}

\newcommand{\BBE}{{\mathbb{E}}}

\newcommand{\sym}{\mathop{\rm sym}\nolimits}
\newcommand{\const}{\mathop{\rm const}\nolimits}
\newcommand{\gh}{\mathfrak{h}}
\newcommand{\cE}{\mathcal{E}}
\newcommand{\cEmS}{{\cE}_{{\min}\varSigma}}

\newcommand{\vph}{\varphi}
\newcommand{\iL}{\mathrm{L}}

\DeclareMathOperator{\ad}{ad}
\newcommand{\arctg}{\arctan}

\newcommand{\by}[1]{\textrm{{#1}}}
\newcommand{\jour}[1]{\textrm{{#1}}}
\newcommand{\vol}[1]{\textrm{{#1}}}
\newcommand{\book}[1]{\textrm{{#1}}}

\begin{document}

\title[Minimal surfaces assigned to contact symmetries]%
{Minimal surfaces associated with nonpolynomial contact symmetries}
\date{April 6, 2006}

\author{Arthemy~V.~Kiselev}
\thanks{\textit{Permanent address}:
  %%%\address{
Department of Higher Mathematics,\ Ivanovo State
Power University,\ Rabfa\-kov\-ska\-ya str.\ 34, Ivanovo 153003, Russia.}
\thanks{\textit{Current address}:
  %%%\curraddr{
Department of Physics, Middle East Technical University,
06531 Ankara, Turkey.
\textit{E-mail}\textup{:
   %%%\email{
\texttt{arthemy\symbol{"40}newton.physics.metu.edu.tr}}}
\thanks{\jour{Fundam.\ Appl.\ Math.} \vol{12} (2006) no.3-4
`Hamiltonian \& Lagrangian systems and Lie algebras,' P.~75--82.}
\thanks{Partially supported by University of Lecce grant no.650~CP/D}

%\address{Department of Higher Mathematics,\ Ivanovo State
%Power University,\ 34 Rabfakovskaya str., 153003 Ivanovo, Russia}
%\curraddr{Department of Physics, Middle East Technical University,
%06531 Ankara, Turkey.}
%\email{arthemy@newton.physics.metu.edu.tr}

  % Symmetries and conservation laws

\begin{abstract}
Two infinite sequences of minimal surfaces in space are constructed using
symmetry analysis. In particular, explicit formulas are obtained for
the self\/-\/intersecting minimal surface that fills the trefoil knot.
\end{abstract}

\maketitle
\noindent{UDC 514.763.85, 517.972.6}

\subsection*{Introduction}
In this paper we consider the Euler\/--\/Lagrange
minimal surface equation
\begin{equation}\label{EMS}
\cEmS=\bigl\{(1+u_y^2)\,u_{xx}-2u_xu_yu_{xy}+(1+u_x^2)\,u_{yy}=0\bigr\}
\end{equation}
whose solutions describe two\/-\/dimensional minimal surfaces
$\varSigma\subset\mathbb{E}^3$
in nonparametric form~$\varSigma=\{ z=u(x, y)\}$,
here $x$,\ $y$,\ $z$ are the Cartesian coordinates.
We construct two infinite sequences of the minimal surfaces
related to nonpolynomial contact symmetries of Eq.~\eqref{EMS}.

\begin{rem}
\textup{Although the graphs of solutions for Eq.~\eqref{EMS} determine
the minimal surfaces only locally such that the projections
of their tangent planes to~$0xy$ are nondegenerate, this is not
restrictive for our reasonings. The minimal surfaces constructed
in section~\ref{SecSurf} are self\/-\/intersecting, being in fact
described by multi\/-\/valued solutions of Eq.~\eqref{EMS} and
admitting singular points.}
\end{rem}

The paper is organized as follows.
In section~\ref{SecSym} we describe the generators and the commutator
relations of the contact symmetry algebra for~$\cEmS$.
We provide examples of the contact non\/-\/point generators and
indicate the recursion operators for the commutative Lie subalgebra
of~$\sym\cEmS$.
In section~\ref{SecSurf} we show that any surface which is invariant
w.r.t.\ a contact non\/-\/point symmetry flow is always a plane,
although non\/-\/planar minimal surfaces in space are assigned to the
same symmetry generators by the inverse Legendre transformation.
Thus we construct two sequences of the minimal surfaces
associated with nonpolynomial contact symmetries of~$\cEmS$;
one of the sequences starts with the
helicoid~\cite{Nitsche} in~$\BBE^3$.
The recursions for the symmetries provide discrete
transformations between the surfaces, while the generators
themselves determine their continuous transformations.
In particular, we obtain explicit formulas for the self\/-\/intersecting
minimal surface~$\varSigma_6$ that fills the trefoil knot;
this surface succeeds the helicoid~$\varSigma_5$ with respect to
the recursion relations.

\begin{rem}
The very idea to construct an infinite sequence of minimal surfaces
with no restrictions upon the boundary conditions can be easily
fulfilled by using the Ennepert\/--\/Weierstrass
representation~\cite{Nitsche} that assigns the surfaces to arbitrary
complex\/-\/analytic functions. The objective of this note is that the
geometry of Eq.~\eqref{EMS} suggests a \emph{natural}
discrete proliferation scheme based on the symmetry approach.
%Also, we exploit an important feature of the minimal surface equation,
%which is polynomial but admits a non\/-\/polynomial irrational source
%symmetry~$\vph_5$ (see Example~\ref{ExampleStart}) that starts the two
%sequences of surfaces.

The resulting surfaces contained in Appendix~\ref{Append} seem to be
relevant in Natural sciences (hydromechanics, bionics, or chemistry);
strangely, these particular solutions given in parametric
representation are not met in classical textbooks and reviews on the
topic~\cite{Nitsche, Osserman}.
\end{rem}

\section{Contact symmetries of the minimal surface
equation}\label{SecSym}
The Ennepert\/--\/Weierstrass representation~\cite{Nitsche}
yields that the symmetry group of the minimal surface equation
is the conformal group, which is a semi\/-\/direct product of
the M\"obius group and the component that corresponds to
complex\/-\/analytic functions (also, the full
symmetry group incorporates
the dilatation). The M\"obius subgroup corresponds to the group of
rotations of~$\BBE^3$ owing to the isomorphism
$\mathfrak{sl}(2)\simeq\mathfrak{so}(3)$. In this section we interpret
the above assertion in view of the Legendre transformation that
brings Eq.~\eqref{EMS} to linear form.

  %\begin{theor}[\cite{Legendre}]\label{LegendreSt}
Let us recall that
equation~\eqref{EMS} is mapped to the linear elliptic equation
\begin{equation}\label{eqGianni}
\mathfrak{L}(\cEmS)=\bigl\{
(1+p^2)\,\phi_{pp}+2pq\,\phi_{pq}+(1+q^2)\,\phi_{qq}=0\bigr\}
\end{equation}
by the Legendre transformation
\[  %$   %  \begin{equation}\label{LegT}
\mathfrak{L}=\{\phi=xu_x+yu_y-u,~p=u_x,~q=u_y\}.
\]  %$   %\end{equation}
The inverse Legendre transformation
$\mathfrak{L}^{-1}=\{ x=\phi_p$, $y=\phi_q$, $u=p\phi_p+q\phi_q-\phi\}$
assigns the minimal surfaces~$\varSigma$ in parametric form
to solutions of Eq.~\eqref{eqGianni}.
   %\end{theor}

Each symmetry of Eq.~\eqref{eqGianni} corresponds to a symmetry
transformation of Eq.~\eqref{EMS}.
Recall that the determining relation
$\iL_\vph(F)=0$ on $\cE=\{F=0\}$
   %$\smash{\bar\ell_F(\varphi)}=0$
for the infinitesimal symmetries~$\varphi$
of any linear differential equation~$\cE$ coincides
with the equation itself, here $\iL_\vph$ is the evolutionary
vector field with the generator~$\vph$ (see~\cite{Olver}).
Therefore it is quite natural that the symmetry algebra of
Eq.~\eqref{EMS} incorporates the set of solutions
\begin{equation}\label{Identif}
\vph(u_x, u_y)=\phi(p, q)
\end{equation}
of the linear equation~\eqref{eqGianni}.
Hence follows the description of contact symmetry algebra
for the minimal surface equation~\eqref{EMS}.

\begin{prop}\label{AlgSymEMS}
The Lie algebra $\sym\cEmS$ of contact symmetries of
the minimal surface equation~\eqref{EMS}
is generated by solutions $\varphi(u_x$,\ $u_y)$
of Eq.~\eqref{eqGianni}, in particular, by the shift
$\varphi_1=1$ and the translations $\varphi_2^i=u_{x^i}$
along $x^1\equiv x$ and $x^2\equiv y$, by the rotations
$\varphi_3^{12}=yu_x-xu_y$ and $\varphi_3^i=x^i+uu_{x^i}$,
here $i=1$\textup{,}\,$2$, and by the
dilatation~$\varphi_4=u-xu_x-yu_y$.
\end{prop}

  %We note that Eq.~\eqref{eqGianni} and its solutions except the point
  %symmetries $\vph_1$ and $\vph_2^i$ were missed in~\cite{Bila}.

\begin{exmp}\label{ExampleStart}
In~\cite{Prague2004} two infinite sequences of the
symmetry generators~$\vph(u_x$,\ $u_y)$ for
the minimal surface equation were constructed.
It was postulated that the functions~$\vph$
are polynomial in~$u_y$; then for each degree~$k\geq0$ of the polynomials
there are two solutions. The initial terms of these sequences are
\begin{align*}
\vph_1&=1, \quad \vph_2^1=u_x, &
\vph_2^2&=u_y, \quad \vph_5=u_y\arctg u_x,\\
\vph_6&=\frac{u_xu_y^2}{1+u_x^2}+\arctg u_x, &
  \vph_{7}&=\frac{u_y^2}{1+u_x^2}-u_x\arctg u_x,\\
\vph_{8}&=\frac{u_xu_y^3}{(1+u_x^2)^2}+\frac{3}{2}\cdot\frac{u_xu_y}{1+u_x^2}, &
  \vph_{9}&=\frac{u_x^2-1}{(1+u_x^2)^2}\cdot
  u_y^3-\frac{3u_y}{1+u_x^2}.
\end{align*}
The generators $\vph_k$ depend rationally on~$u_x$ for all~$k\geq8$.
We conjecture that none of the contact symmetries~$\vph_k$
is Noether whenever~$k\geq5$.
\end{exmp}

The identification~\eqref{Identif} yields that a sequence of
solutions~$\phi(p,q)$ is obtained whenever a recursion for the contact
non\/-\/point symmetries~$\vph(u_x,u_y)$ of Eq.~\eqref{EMS}
is known. We claim that three local recursion operators
for this component of~$\sym\cEmS$ are determined by
the adjoint representation of the symmetry algebra itself.
Now we study these aspects in more detail.

The commutation relations for the seven point symmetries
$\varphi_1$, $\ldots$, $\varphi_4$, see Proposition~\ref{AlgSymEMS},
were derived in~\cite{Bila}.
Let us indicate the commutation properties of the contact symmetries
that originate from Eq.~\eqref{eqGianni}.

\begin{lem}
Assume that $\vph'(u_x, u_y)$ and $\vph''(u_x, u_y)$ are the generators
of evolutionary vector fields $\iL_{\vph'}$
and~$\iL_{\vph''}$. Then their Jacobi bracket $\{\vph'$\textup{,}
$\vph''\}$ is always trivial.
\end{lem}

%\noindent\textit{Proof.}\quad
%In local coordinates, the Jacobi bracket $\{\vph'$\textup{,}
%$\vph''\}=\iL_{\vph'}(\vph'')-\iL_{\vph''}(\vph')$ is
%\begin{multline*}
%\{\vph',\vph''\}=
% \frac{\dd\vph'}{\dd u_x}u_{xx}\cdot\frac{\dd\vph''}{\dd u_x}+
% \frac{\dd\vph'}{\dd u_y}u_{xy}\cdot\frac{\dd\vph''}{\dd u_x}+
% \frac{\dd\vph'}{\dd u_x}u_{xy}\cdot\frac{\dd\vph''}{\dd u_y}+
% \frac{\dd\vph'}{\dd u_y}u_{yy}\cdot\frac{\dd\vph''}{\dd u_y} -
%\text{v.\,v.} = 0.
%\end{multline*}
%Recall that the commutator of two fields is well defined.
%This completes the proof.\hfill$\Box$

\begin{prop}\label{CommutatorSt}
%\begin{enumerate}
%\item
All contact symmetries $\vph(u_x$, $u_y)\in\sym\cEmS$ of the minimal
surface equation~$\cEmS$ commute.
%\item
  %Suppose further that $\varphi(u_x$, $u_y)\in\sym\cEmS$ is a
  %symmetry\textup{;} then
Also, the following relations hold\textup{:}
\begin{align*}
\{\varphi_3^{12},\varphi\} &= u_x\,\frac{\partial\varphi}{\partial u_y}
- u_y\,\frac{\partial\varphi}{\partial u_x}, \qquad
\{\varphi_4,\varphi\} = -\varphi,
\\
\{\varphi_3^{i},\varphi\} &= -u_{x^i}\varphi +
  (1+u_{x^i}^2)\,\frac{\partial\varphi}{\partial u_{x^i}} +
  u_xu_y\,\frac{\partial\varphi}{\partial u_{x^{3-i}}}.
\end{align*}
The Lie subalgebra~$\gh$ generated by
the solutions $\varphi(u_x$, $u_y)$
of Eq.~\eqref{eqGianni} is the radical of the contact symmetry
algebra~$\sym\cEmS$ for Eq.~\eqref{EMS}.
%\end{enumerate}
\end{prop}

\begin{cor}\label{Recursions}
The mappings $\ad_{\vph_3^{12}}$ and $\ad_{\vph_3^{i}}\colon\gh\to\gh$
define the local recursion operators on
the Lie subalgebra~$\gh\subset\sym\cEmS$.
\end{cor}

\begin{rem}\label{RemDiagr}
The symmetries introduced in Example~\ref{ExampleStart} are
proliferated according to the following diagram~\cite{YS2005}:
\[
\begin{CD}
\vph_5 @>{\ad_{\vph_3^2}}>> \vph_6 @>{\ad_{\vph_3^1}}>> \vph_7
 @>{-\frac{1}{2}\ad_{\vph_3^{12}}}>> \vph_8-\frac{1}{2}\vph_5 \\
@V{\ad_{\vph_3^1}}VV @. @. @VV{-\ad_{\vph_3^1}}V \\
\vph_2^2 @>{\ad_{\vph_3^2}}>> \vph_1
@>{-\ad_{\vph_3^1}}>> \vph_2^1 @. \vph_9+\frac{7}{2}\vph_2^2.
\end{CD}
\]
Two infinite sequences of the contact symmetries are
obtained~\cite{Prague2004} from~$\vph_6$
and~$\vph_7$ by multiple application of the recursion~$\ad_{\vph_3^{12}}$.
From the above diagram it follows that
the symmetry~$\vph_5$ is the `seed' generator for both sequences.
\end{rem}

\section{The minimal surfaces associated with the contact
       symmetries}\label{SecSurf}
Continuing the line of reasonings, we see
that any contact symmetry $\vph(u_x,u_y)\in\gh$ of Eq.~\eqref{EMS}
determines the minimal surfaces~$\varSigma$ using two different methods:
\begin{enumerate}
\item
Recall that $\phi(p, q)$ defined in~\eqref{Identif}
is a solution of Eq.~\eqref{eqGianni},
hence the inverse image $\mathfrak{L}^{-1}(\phi)$ of the graph of~$\phi$
with respect to the Legendre
transformation is a minimal surface in parametric representation.
\item
The symmetry reduction $\cEmS\cap\{\vph=0\}$ of Eq.~\eqref{EMS}
by a generator $\vph\in\gh\subset\sym\cEmS$
determines the $\vph$-\/invariant surface in~$\BBE^3$.
\end{enumerate}

\begin{exmp}
Consider the solution $\phi_5=q\arctg p$ of Eq.~\eqref{eqGianni}.
Using the first method, we obtain the helicoid
$\{z=x\tan y\}\subset\BBE^3$ whose axis is~$0y$.
The minimal surface which is invariant w.r.t.\ the
symmetry~$\vph_5=u_y\,\arctan u_x$ is a plane.
This is a particular case of the following general property of
the surfaces.
\end{exmp}

\begin{prop}
Suppose that a minimal surface $\varSigma$ is invariant w.r.t.\ a
contact symmetry $\vph(u_x$, $u_y)$. Then $\varSigma$ is a plane.
\end{prop}

\begin{proof}
Consider the constraint $\vph=0$. By the implicit function theorem, we
have $u_y=f(u_x)$ almost everywhere. Therefore,
$u_{yy}=\left(f'(u_x)\right)^2\cdot u_{xx}$ and from
Eq.~\eqref{EMS} we obtain the equation
\[
\Bigl(
(1+u_x^2)\cdot\left(f'(u_x)\right)^2 -
2u_x\,f(u_x)\cdot f'(u_x) + (1+f^2(u_x))\Bigr) \cdot u_{xx} = 0.
\]
Hence either $u_{xx}=0$ and we have $u_x=\const$, $u_y=f(u_x)=\const$,
or $u_x$ satisfies the algebraic equation whose solutions
are $u_x=\const$ and therefore $u_y=f(u_x)=\const$ again.
\end{proof}

In what follows we focus on the first method for constructing the minimal
surfaces, that is, $\varSigma=\mathfrak{L}^{-1}(\phi(p,q))$.
Using Proposition~\ref{CommutatorSt} and Corollary~\ref{Recursions},
we obtain two sequences of solutions~$\phi_k(p,q)$
of Eq.~\eqref{eqGianni}, which are polynomial in~$q$ and which
are nonpolynomial in~$p$ for~$k\geq5$.
By Remark~\ref{RemDiagr}, both sequences are obtained from the
generating section~$\phi_5(p,q)$.
First let us list the solutions $\phi_{10}$, $\phi_{11}$
and $\phi_{12}$, $\phi_{13}$ of Eq.~\eqref{eqGianni},
which are polynomials in~$q$ of degrees~$4$ and~$5$, respectively.
We have
\begin{align*}
\phi_{10}&=
  \frac{p^3-3p}{{(1+p^2)}^3}\cdot q^4
  +\frac{3}{2}\cdot\frac{p^5-2p^3-3p}{{(1+p^2)}^3}\cdot q^2
  -\frac{3}{2}\cdot\frac{p}{1+p^2},\\
\phi_{11}&=
  \frac{3p^2-1}{{(1+p^2)}^3}\cdot q^4
  +\frac{3}{2}\cdot\frac{3p^4+2p^2-1}{{(1+p^2)}^3}\cdot q^2
  +\frac{3}{2}\cdot\frac{p^2}{1+p^2},\\
\phi_{12}&=
  \frac{p^4-6p^2+6}{{(1+p^2)}^4}\cdot q^5
  +\frac{11p^6-49p^4-51p^2+9}{6{(1+p^2)}^4}\cdot q^3
  +\frac{2p^8-3p^6-11p^4-5p^2+1}{2{(1+p^2)}^4}\cdot q,\\
\phi_{13}&=
  \frac{p^3-p}{{(1+p^2)}^4}\cdot q^5
  +\frac{21p^5+2p^3-19p}{12{(1+p^2)}^4}\cdot q^3
  +\frac{3p^7+4p^5-p^3-2p}{4{(1+p^2)}^4}\cdot q.
\end{align*}

The resulting surfaces can be easily visualized using standard software,
e.g., by modifying this sample code for \textsc{Maple}:
\begin{quote}
\begin{verbatim}
phi6:=p*q^2/(1+p^2)+arctan(p);
x:=diff(phi6,p); y:=diff(phi6,q); z:=p*x+q*y-phi6;
plot3d([x,y,z],p=-2..2,q=-2..2,grid=[50,50]);
\end{verbatim}
\end{quote}
We discover that the boundary of the surface
$\varSigma_6=\mathfrak{L}^{-1}(\phi_6)$ is supported by the trefoil
knot such that the self\/-\/intersecting surface~$\varSigma_6$
fills its interiour. The regular minimal surface
$\varSigma_7=\mathfrak{L}^{-1}(\phi_7)$ resembles a flying bird.
The surfaces~$\varSigma_8$ and~$\varSigma_{13}$ are propeller\/-\/like.
The family~$\varSigma_{9-12}$ provides the shapes of pearl shells;
the origin is a singular point for them, and their self\/-\/intersections
divide the space~$\BBE^3$ in different number of cells for different
subscripts~$k$. The corresponding pictures are given in
Appendix~\ref{Append}.

%\begin{rem}
%The generators of the Lie algebra $\sym\cEmS$,
%which were described in Proposition~\ref{AlgSymEMS},
%provide the continuous symmetry transformations~$u_t=\vph$ of
%the surfaces $\varSigma=\{z=u(x$,\ $y)\}$.
%The recursion operators, see Corollary~\ref{Recursions},
%define the discrete transformations between the surfaces~$\varSigma$.
  %Hence we conclude that each surface we have defined belongs
  %to an infinite\/-\/parametric continuous family of minimal surfaces
  %in space.
%\end{rem}

\subsection*{Acknowledgements}
The author thanks A.\,Klyachko, D.\,Pelinovskiy, and
R.\,Vitolo for helpful discussions.
This research was partially supported by University of Lecce
grant no.650~CP/D.
A part of this research was carried out while the author was
visiting at University of Lecce.

   %%%\rightline{Translated by \textsc{the Author}.}

\newpage

\appendix
\section{The plots of minimal surfaces}\label{Append}

\subsection{The surface $\varSigma_6$: trefoil knot}
\[
\includegraphics{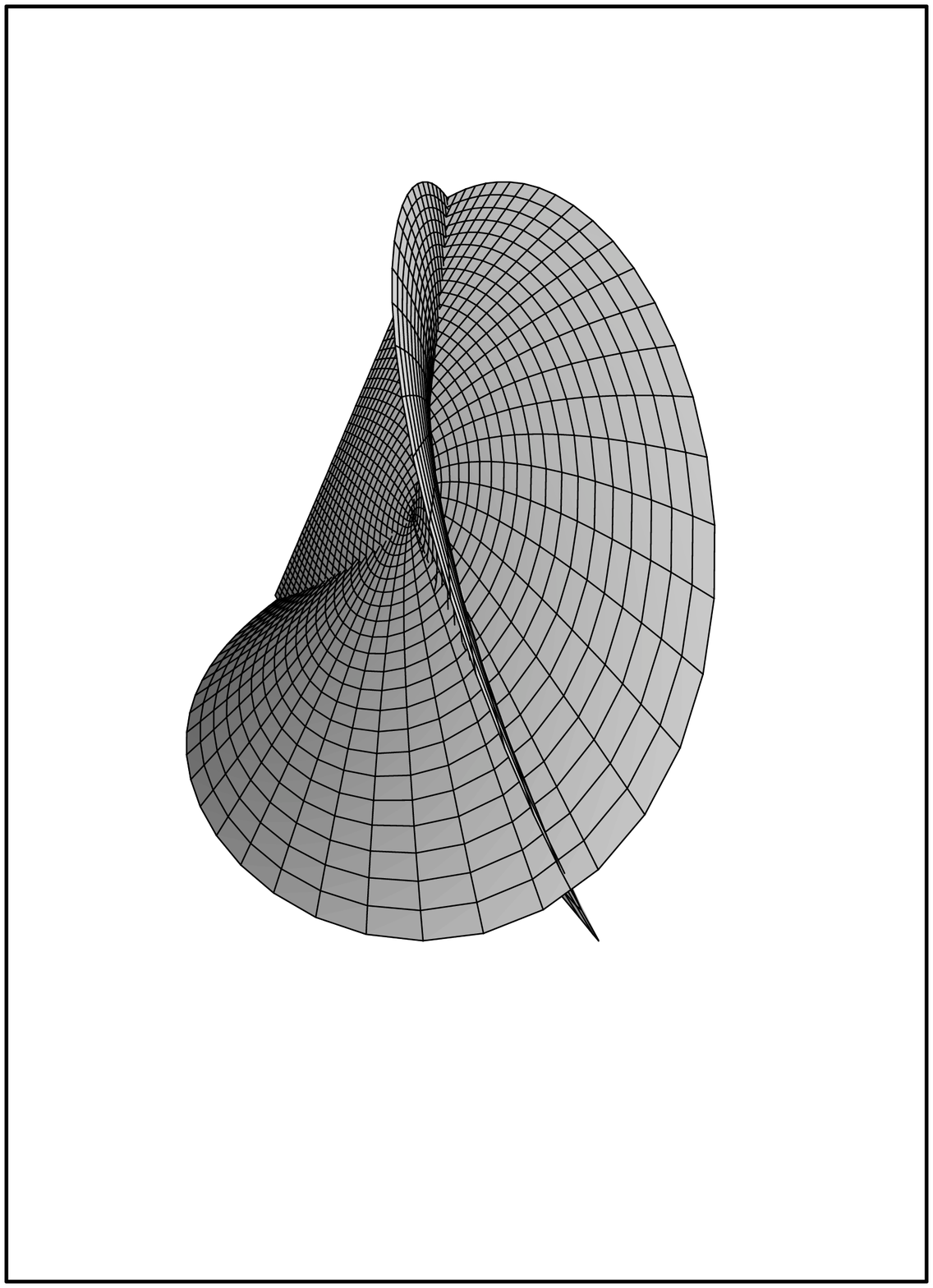}
\]

\subsection{The surface $\varSigma_7$: flying bird}
\[
\includegraphics{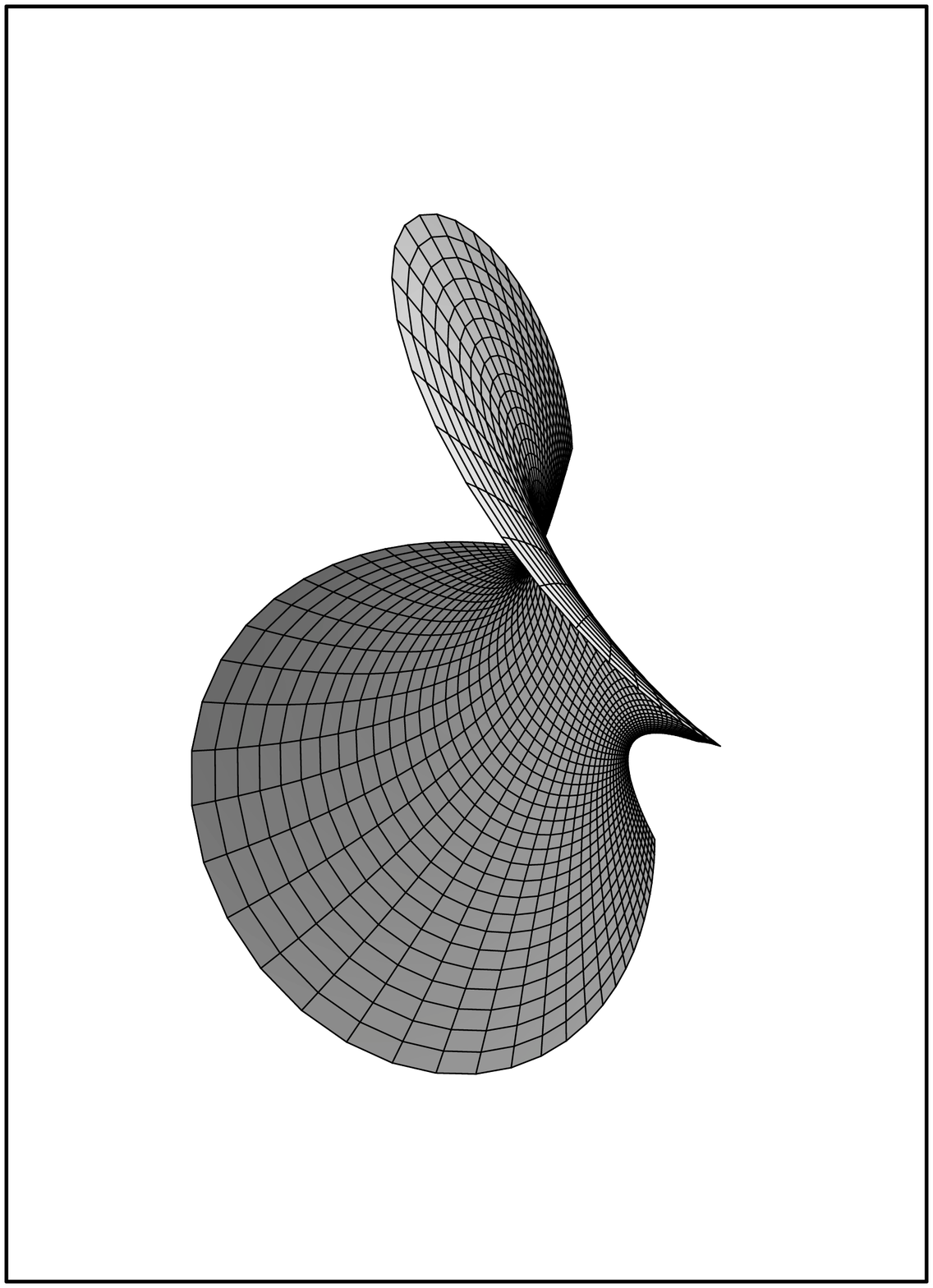}
\]

\subsection{The surfaces $\varSigma_8$ and $\varSigma_{13}$: propellers}
\[
\includegraphics{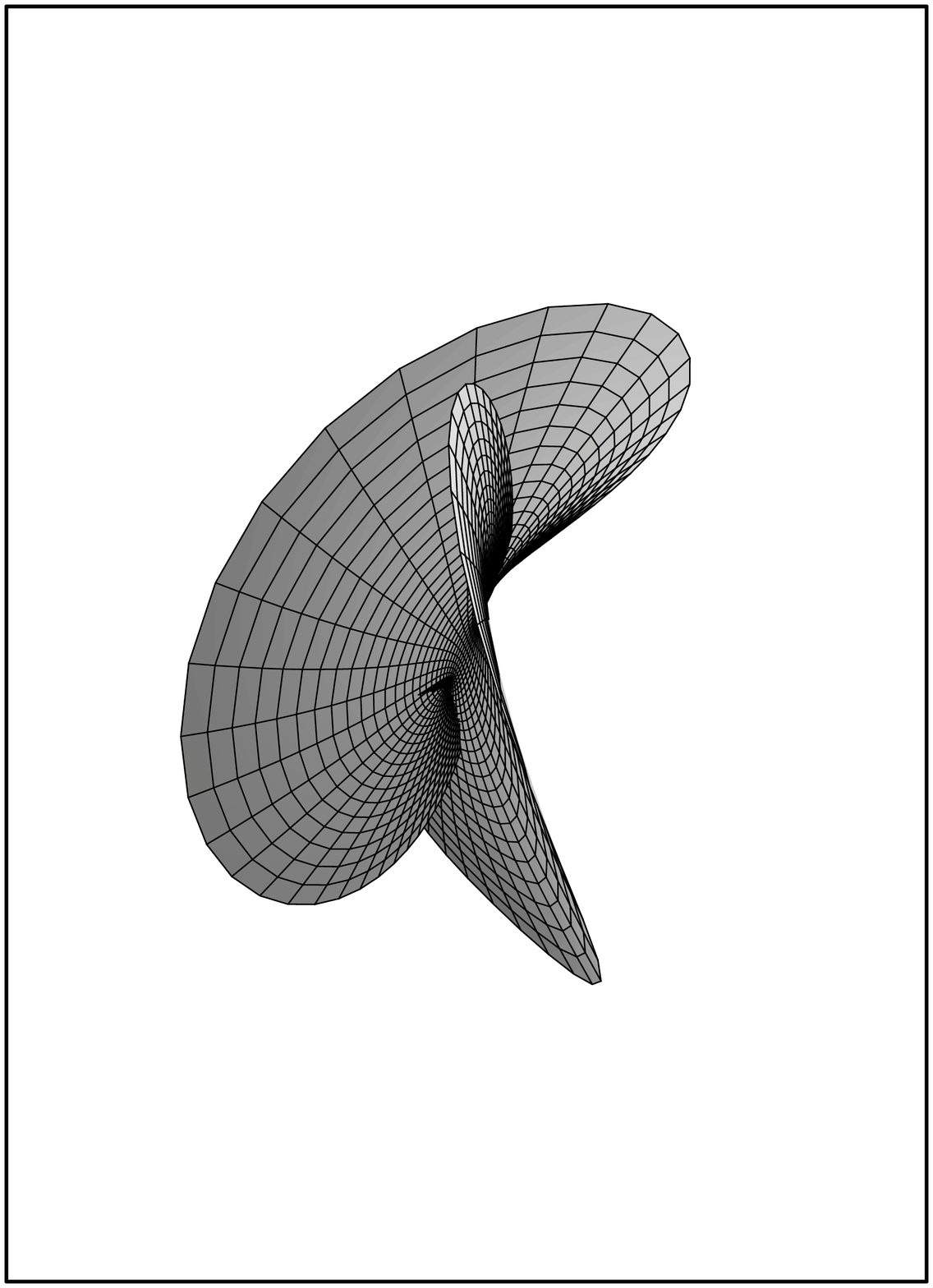}
\]
The surface~$\varSigma_{13}$:
\[
\includegraphics{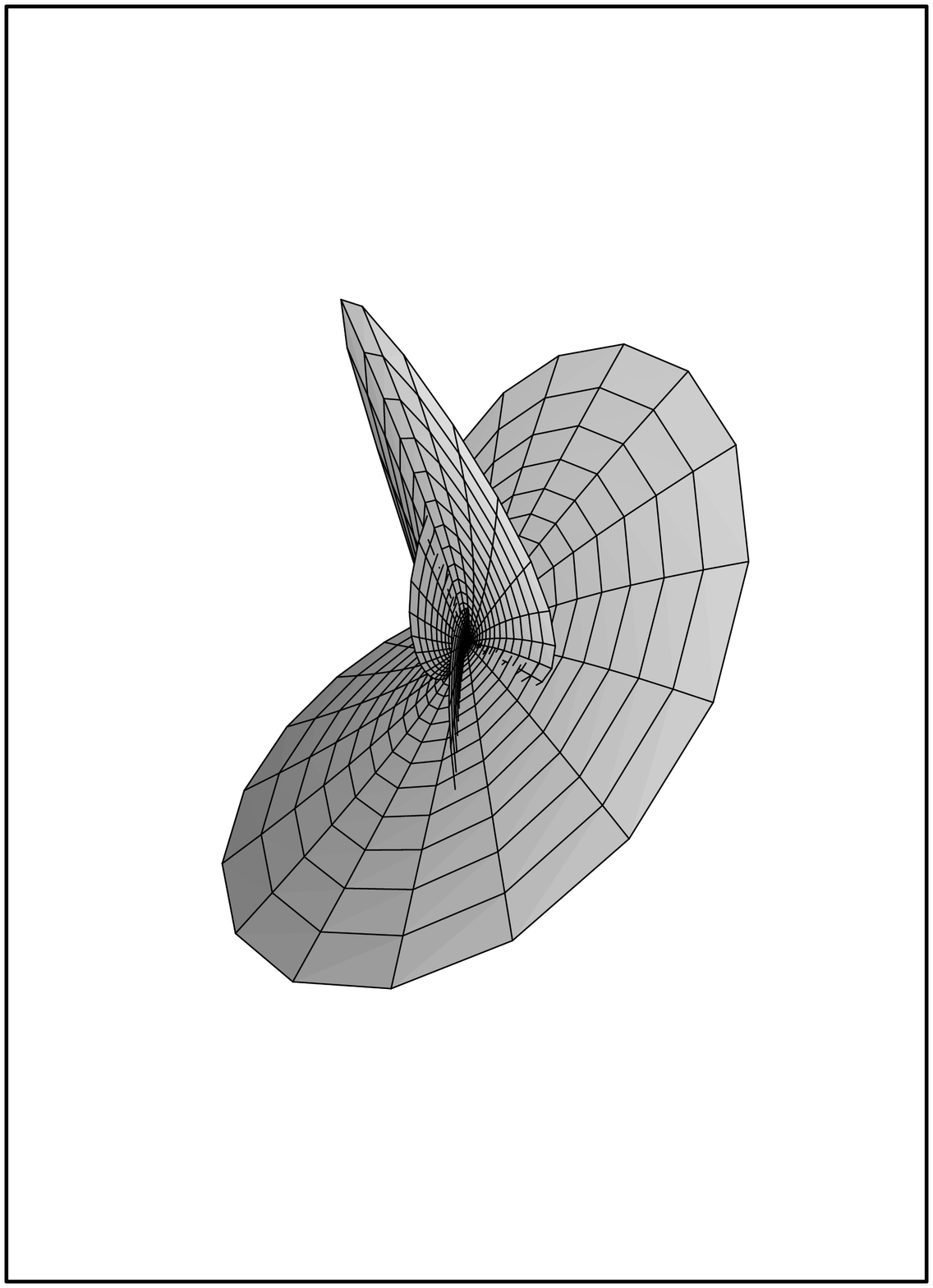}
\]

\subsection{The surfaces $\varSigma_{9-12}$: pearl shells}
\[
\includegraphics{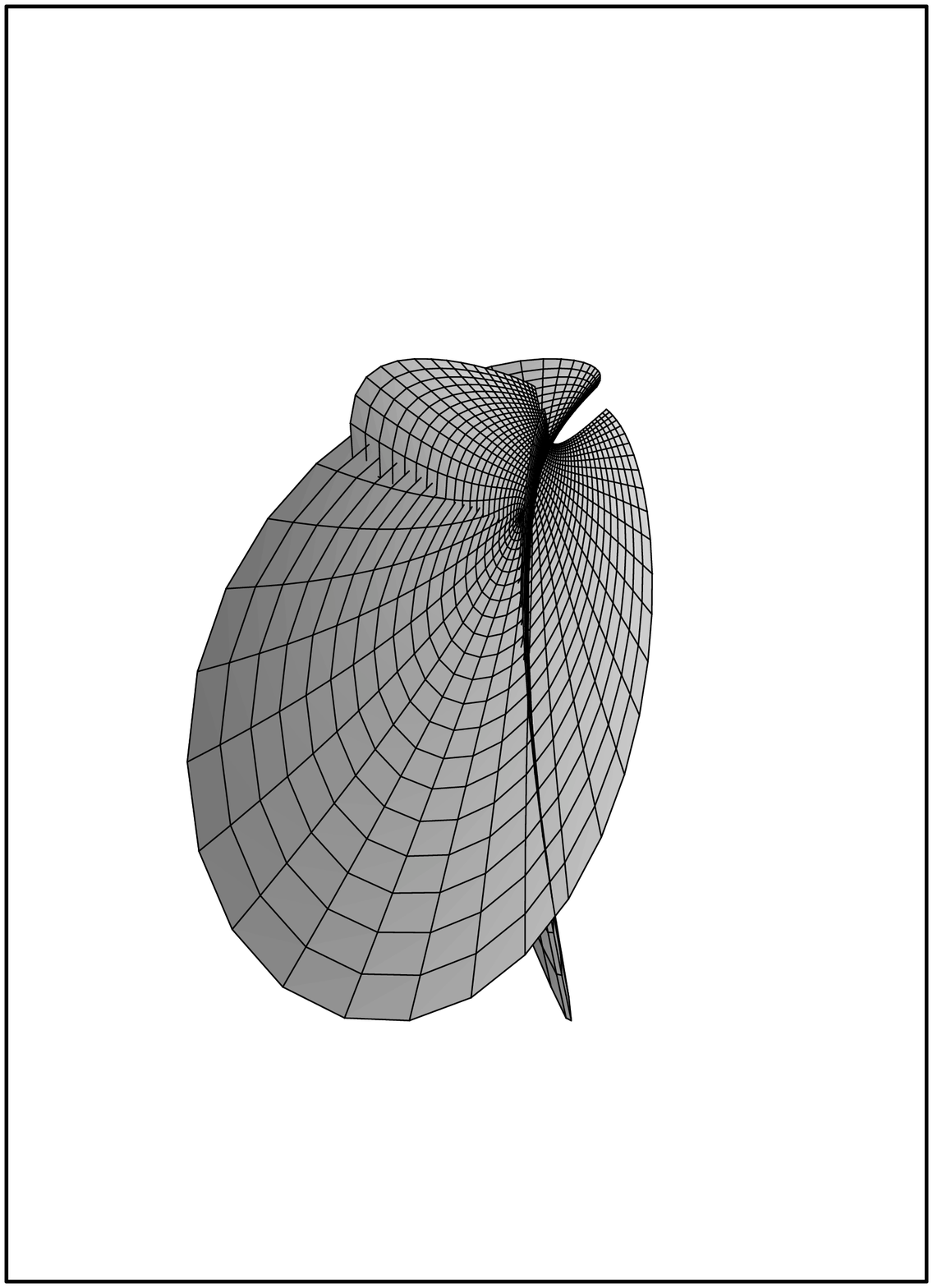}
\]
The surface $\varSigma_{10}$:
\[
\includegraphics{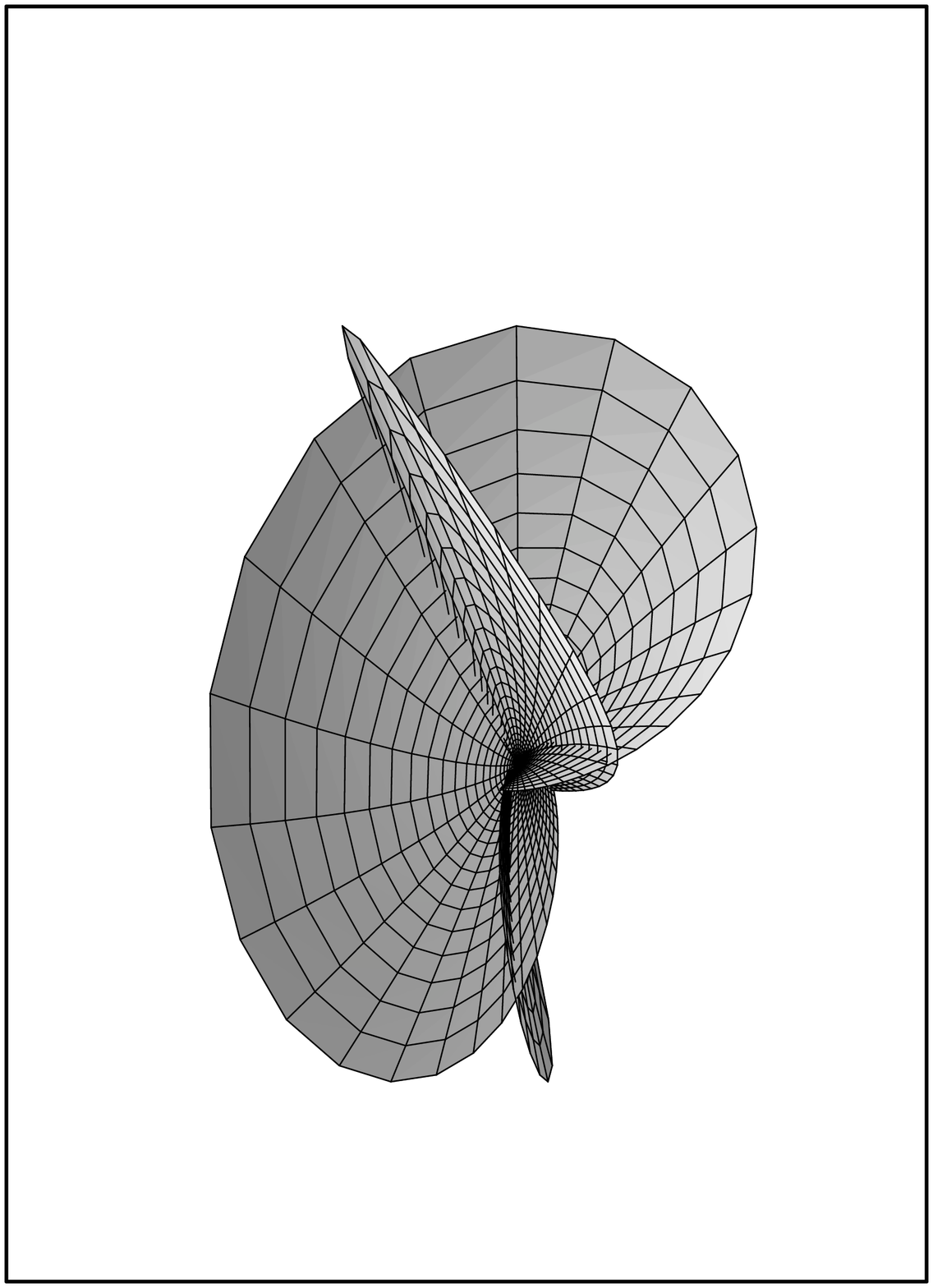}
\]
The surface $\varSigma_{11}$:
\[
\includegraphics{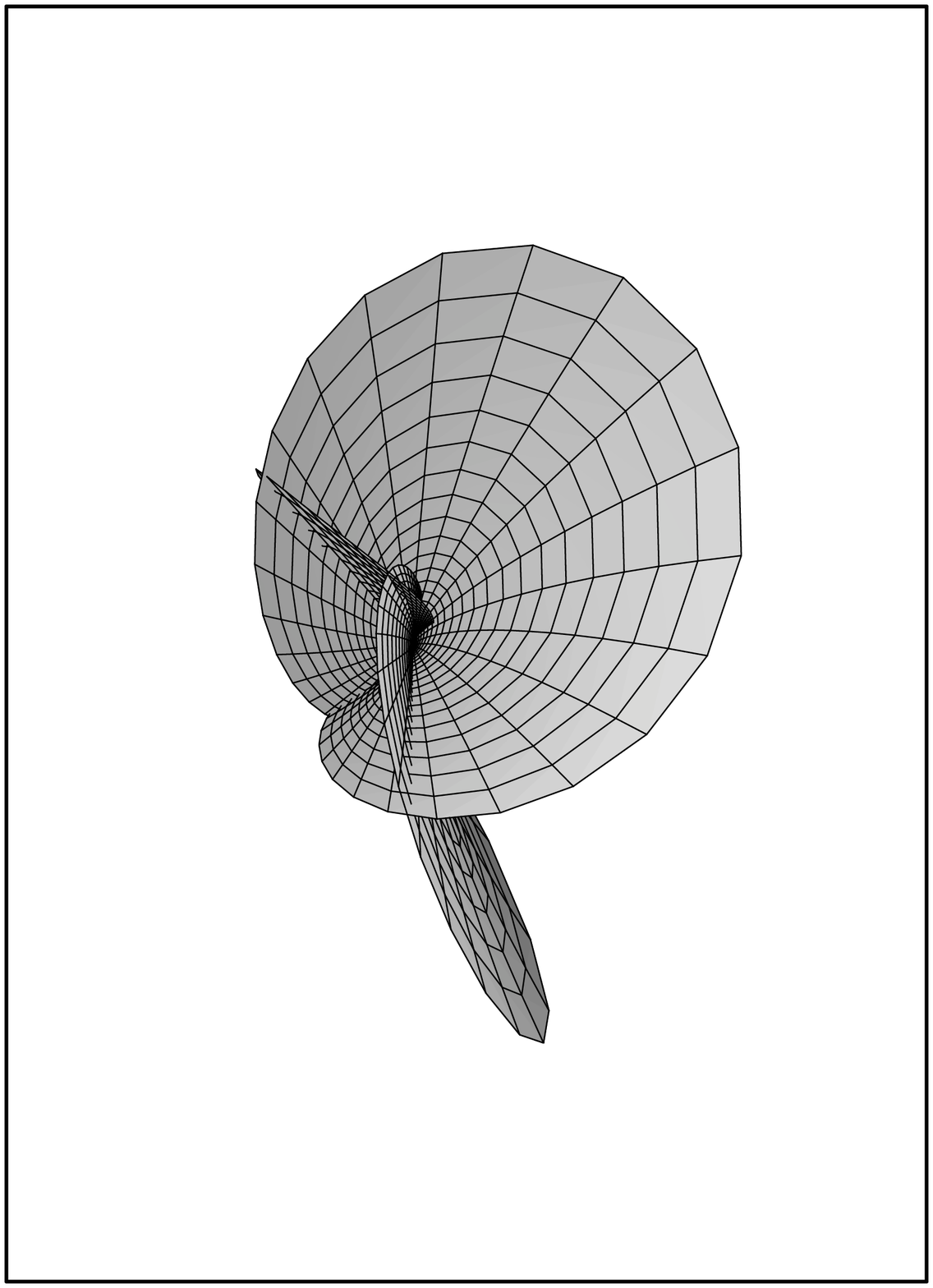}
\]
The surface $\varSigma_{12}$:
\[
\includegraphics{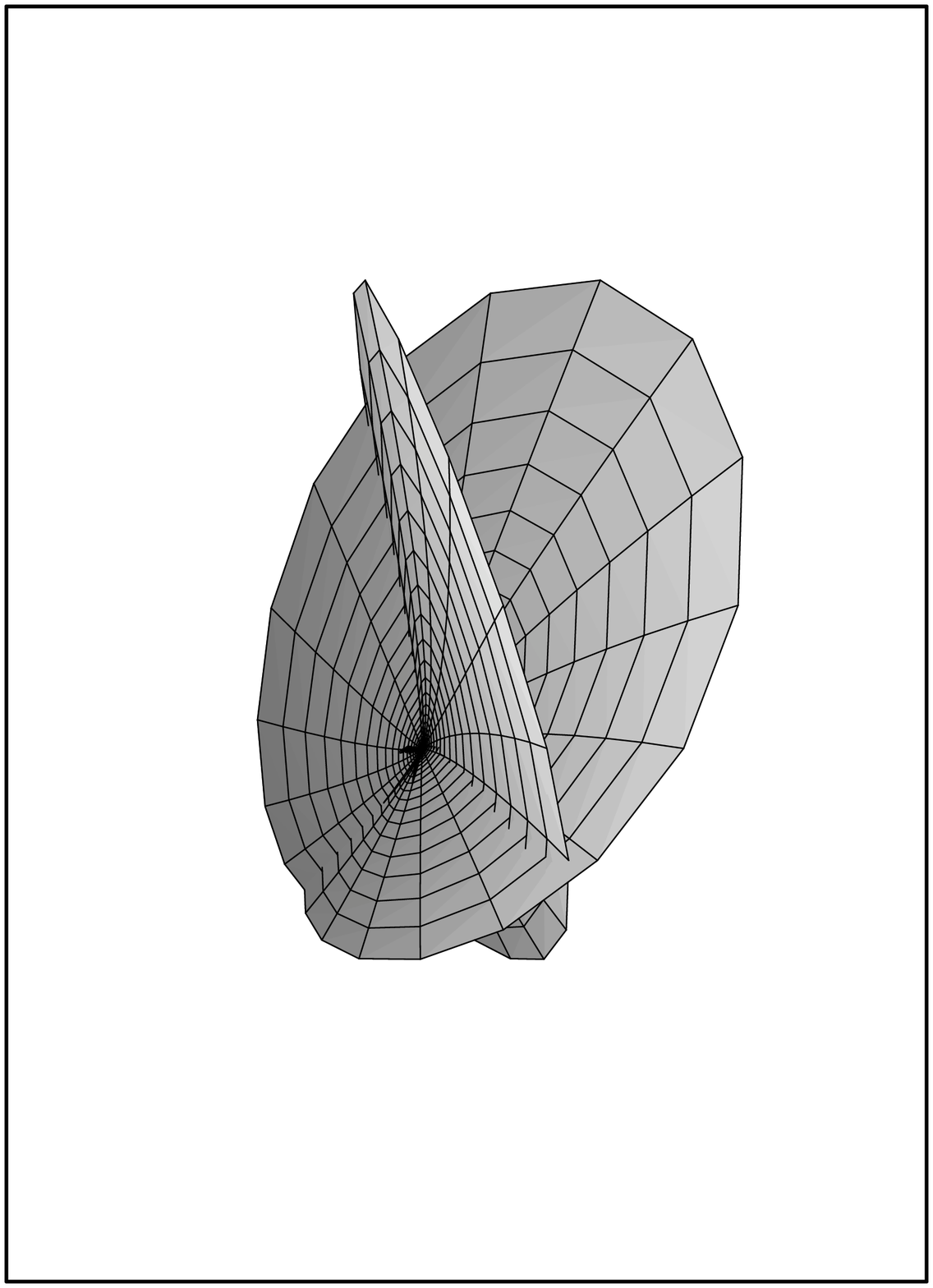}
\]

\begin{thebibliography}{9}\normalsize

\bibitem{Bila}
\by{N.~B\^{\i}l\u{a}}, Lie groups applications to minimal surfaces PDE,
\jour{Diff.\ Geometry~-- Dynam.\ Systems} \vol{1} (1999) no.1, 1--9.

%\bibitem{ClassSym}
%\by{A.V. Bocharov, V.N. Chetverikov, S.V. Duzhin et al.},
%\book{Symmetries and Conservation Laws for Differential Equations of
%Mathematical Physics}. AMS, Providence, RI, 1999.
%I. Krasil'shchik and A. Vinogradov (eds).

\bibitem{YS2005}
\by{A.V. Kiselev}, On a symmetry reduction of the minimal surface
equation, Proc.\ XXVII Conf.\ of Young Scientists, April 2005.
Faculty of Mathematics and Mechanics, Lomonosov MSU, Moscow (2006),
67-71.

  %\bibitem{NoteNonpolynomial}
  %\by{A.V. Kiselev, G. Manno}, On existence of nonpolynomial structures
  %admitted by polynomial differential equations, \jour{Bulletin of
  %ISPU} (2004) 4, 30--31.

\bibitem{Prague2004}
\by{A.V. Kiselev, G. Manno}, On the symmetry structure of the minimal
    surface equation, Proc.\ IX Conf.\ `Differential Geometry and Its
    Applications,' August 2004. Charles Univ., Prague (2005), 483-490.
    \texttt{arXiv:math.DG/0410557}

%\bibitem{Legendre}
%\by{Legendre~A.}, M\'emoire sur l'int\'egration de quelques
%\'equ\-a\-ti\-ons aux
%dif\-f\'e\-ren\-ces partielles, \jour{M\'em.\ Acad.\
%Roy.\ Sci.\ Paris} (1789), 309--351.

\bibitem{Nitsche}
\by{J.C.C.~Nitsche}, \book{Vorlesungen \"uber Minimalfl\"achen}.
Springer, Berlin, 1974.

\bibitem{Olver}
\by{P.J.~Olver}, \book{Applications of Lie groups to differential
equations}, $2$nd ed.,  %Graduate Texts in Mathematics \vol{107},
Springer, Berlin, 1993.

\bibitem{Osserman}
\by{Osserman R.}, \book{A survey of minimal surfaces}.
Van Nostrand Reinhold, NY, 1969.

\end{thebibliography}
\end{document}